# Mixed finite element and TPSA finite volume methods for linearized elasticity and Cosserat materials

J. M. Nordbotten[1,2], W. M. Boon[2], O. Duran[1], E. Keilegavlen[1]

1: Center for Modeling of Coupled Subsurface Dynamics,
Department of Mathematics, University of Bergen
Norway

2: Norwegian Research Center (NORCE)
Bergen, Norway

The equations of linearized elasticity are most commonly discretized by finite element methods, however in the context of subsurface applications there are several advantages to considering finite volume methods. Such advantages include compatibility with the grid structure used for flow equations and explicit and robust representation of the traction at internal surfaces, including fractures.

In previous work, we have developed extensions of the classical multi-point flux approximation (MPFA) methods for flow to multi-point stress approximation (MPSA) methods for elasticity. A consistent finite volume method for linearized elasticity with a simpler, two-point stress approximation (TPSA) has only recently been discovered, inspired by linear Cosserat materials.

Cosserat theory of elasticity is a generalization of classical elasticity that allows for asymmetry in the stress tensor by taking into account micropolar rotations in the medium. The equations involve a rotation field and associated "couple stress" as variables, in addition to the conventional displacement and Cauchy stress fields.

In recent work, we derived a mixed finite element method (MFEM) for the linear Cosserat equations that converges optimally in these four variables. The drawback of this method is that it retains the stresses as unknowns, and therefore leads to relatively large saddle point system that are computationally demanding to solve.

As an alternative, we developed a finite volume method in which the stress variables are approximated using a minimal, two-point stencil. The system consists of the displacement and rotation variables, with an additional "solid pressure" unknown.

Both the MFEM and TPSA methods are robust in the incompressible limit and in the Cauchy limit, for which the Cosserat equations degenerate to classical linearized elasticity. We report on the construction of the methods, their a priori properties, and compare their numerical performance against an MPSA finite volume method.

**Introduction**

Cosserat theory forms a generalization of classical linearized elasticity by relaxing the symmetry of the stress tensor (Cosserat and Cosserat 1909, Gurtin 1982, Truesdell and Noll 1965). In particular, the asymmetric part of the stress is related to a secondary field that represents local, or micropolar, rotations. In the case of granular media, for example, this field describes rotations of the grains at a small length scale, which follow a different dynamic than the rotation of the medium at the domain scale.

The generalization to Cosserat materials has proven key in designing robust and simple numerical methods for linearized elasticity and Stokes flow. In fact, it is through this perspective that the first consistent two-point stress approximation (TPSA) finite volume method was derived for these systems on general grids (Nordbotten and Keilegavlen 2024). In this work, we assess the performance of this finite volume method and compare it to the recently proposed four-field finite element method for Cosserat materials (Boon, Duran and Nordbotten 2024), which is an extension of the seminal three-field finite element method for elasticity (Arnold, Falk and Winther 2007). For reference, we also include the multi-point stress approximation method (MPSA) (Keilegavlen and Nordbotten, Finite volume methods for elasticity with weak symmetry 2017).

**1. Theory and Methods**

In this section, we present the equations governing linear Cosserat materials, as well as the discretizations considered in this paper. We describe the system for a three-dimensional computational domain $\Omega$ and note the differences with the two-dimensional analogue where needed.

First, letting $\sigma$ denote the Cauchy stress tensor, the balance of linear momentum is stated as:
$$-\nabla \cdot \sigma = f \qquad (1)$$
in which $f$ describes a body force.

Second, given Lamé parameters $\lambda$ and $\mu$, we introduce Hooke's law for homogeneous, isotropic media, given by:
$$\sigma = 2\mu(\nabla u + S^* r_s) + \lambda(\nabla \cdot u)I \qquad (2)$$
Here $u$ is the material displacement, and $r_s$ is the (local) material rotation. The operator $S^*: \mathbb{R}^3 \to \mathbb{R}^{3\times3}$ is given by:
$$S^* r_s := \begin{pmatrix} 0 & -r_3 & r_2 \\ r_3 & 0 & -r_1 \\ -r_2 & r_1 & 0 \end{pmatrix}$$

Third, to measure the asymmetry of a tensor $\sigma$, we define the operator $S: \mathbb{R}^{3\times3} \to \mathbb{R}^3$, for which $S^*$ is the adjoint, such that
$$S(\sigma)_i := \sigma_{i-1,i+1} - \sigma_{i+1,i-1}$$
in which the indices are understood modulo 3. Classically, the Hellinger-Riesner formulation of linearized elasticity is obtained by requiring $S\sigma = 0$. Here we follow the somewhat more general Cosserat description of mechanics, where the *couple stress* $\omega$ is introduced to balance any deviation from a symmetric stress tensor according to the following conservation law:
$$\nabla \cdot \omega = S(\mu^{-1}\sigma) \qquad (3)$$

Fourth, by letting $\ell \geq 0$ denote a *micropolar length scale*, we close the system by relating the couple stress to the local material rotation according to the lowest-order constitutive law:
$$\omega = 2\ell^2 \nabla r_s \qquad (4)$$

The above system is complemented by homogeneous boundary conditions on the displacement and rotation variables:
$$u = 0, \qquad r_s = 0, \qquad \text{on } \partial\Omega. \qquad (5)$$

*Remark* In 2D, the rotational variable $r_s$ is a scalar field since there is only one axis around which rotations are possible. Moreover, the stress $\sigma$ is a $2 \times 2$ tensor field, the asymmetry operator $S: \mathbb{R}^{2\times 2} \to \mathbb{R}$ becomes $S(\sigma) \coloneqq \sigma_{21} - \sigma_{12}$, and its adjoint is given by $S^*r \coloneqq \begin{pmatrix} 0 & -r \\ r & 0 \end{pmatrix}$. Finally, we note that the couple stress $\omega$ is a vector field in 2D, not a tensor field.

We briefly demonstrate that the classical equations for linearized elasticity are recovered from (1)-(4) in the limit $\ell \to 0$. We refer to this as the *Cauchy limit* since it follows the classical arguments by Cauchy regarding the symmetry of the stress tensor. In particular, inserting $\ell = 0$ in (4) implies $\omega = 0$ and, in turn, it follows from (3) that $\sigma$ is symmetric. Substituting (2) in (3) then gives us that $0 = S(\nabla u + S^*r_s)$ and, after a short calculation, we obtain $r_s = -\frac{1}{2}S(\nabla u)$. It follows by direct computation that $\nabla u + S^*r_s = \nabla u - \frac{1}{2}S^*S(\nabla u) = \frac{1}{2}(\nabla u + \nabla u^T)$, which we recognize as the linearized strain $\varepsilon(u)$. Using this identity in (2), we obtain the more conventional form of Hooke's law:
$$\sigma = 2\mu\varepsilon(u) + \lambda(\nabla \cdot u)I. \tag{6}$$

If we additionally consider the incompressible limit $\lambda \to \infty$, then we obtain the equations governing Stokes flow. Paying particular attention to these limits therefore allows us to construct discretization methods that are applicable to a wide class of applications, ranging from Cosserat materials to Stokes fluids. While not the topic of this paper, we also remark that robustness in the limit of $\lambda \to \infty$ is also advantageous when considering poromechanical materials.

Let us return to the general case in which the system (1)-(5) describes linear Cosserat materials in terms of the variables $(\sigma, \omega, u, r_s)$. We will consider and compare two discretization approaches. The first is based on mixed finite elements whereas the second is a finite volume method with minimal stencil. In line with the arguments outlined above, these methods are specifically designed to capture the Cauchy limit $\ell \to 0$ and the incompressible limit $\lambda \to \infty$, albeit by employing different techniques.

## 1.1 A four-field mixed finite element method

We start by considering the lowest-order mixed finite element method proposed in (Boon, Duran and Nordbotten 2024). For this exposition, we will slightly abuse notation by letting $(\cdot,\cdot)_\Omega$ denote the $L^2(\Omega)$-inner product for vector fields $u$ and $\tilde u$ as well as for tensor fields $\sigma$ and $\tilde\sigma$, i.e.
$$(u, \tilde u)_\Omega \coloneqq \int_\Omega u \cdot \tilde u \, dx, \qquad (\sigma, \tilde\sigma)_\Omega \coloneqq \int_\Omega \sigma : \tilde\sigma \, dx.$$

Our aim is to derive a variational formulation of (1)-(4) that leads to a symmetric discretization matrix. First, the momentum balance can directly be tested with test functions $\tilde u$ to yield:
$$(-\nabla \cdot \sigma, \tilde u)_\Omega = (f, \tilde u)_\Omega, \qquad \forall \tilde u. \tag{7}$$
For brevity, we will not specify the function spaces that the solution and test functions belong to in the continuous setting. These details are elaborated in (Boon, Duran og Nordbotten 2024).

Next, we consider Hooke's law (4). Given $d$ the dimension of the domain $\Omega$, let the operator $A: \mathbb{R}^{d\times d} \to \mathbb{R}^{d\times d}$ be such that its inverse is given by $A^{-1}\tau = 2\mu\tau + \lambda(\tau:I)I$. It can be written in closed form as
$$A\sigma \coloneqq \frac{1}{2\mu}\left(\sigma - \frac{\lambda}{d\lambda + 2\mu}(\sigma:I)I\right)$$
This allows us to rewrite (4) as
$$A\sigma = \nabla u + S^*r_s$$
Using integration by parts and the fact that $S^*$ and $S$ are adjoints, we obtain the following variational form of (4):
$$(A\sigma, \tilde\sigma)_\Omega + (u, \nabla \cdot \tilde\sigma)_\Omega - (r_s, S\tilde\sigma)_\Omega = 0, \qquad \forall \tilde\sigma. \tag{8}$$

We now turn our attention to (4). We manipulate this relation such that the gradient on $r_s$ appears as the adjoint of the divergence on the couple stress in equation (3). However, if we simply divide both sides by $\ell^2$, this would introduce a scaling of $\ell^{-2}$ in the system which will make the Cauchy limit difficult to handle. Alternatively, introducing a scaled $\hat{\omega} := \ell^{-2}\omega$ will lead to an inconvenient scaling in (3). To remedy this, we instead define the scaled couple stress as $\omega_\ell := \mu\ell^{-1}\omega$ and divide (4) by $2\mu\ell$. After integration by parts, the weak formulation of (3) then becomes:

$$((2\mu)^{-1}\omega_\ell, \tilde{\omega})_\Omega + (r_s, \nabla \cdot (\ell\tilde{\omega}))_\Omega = 0, \qquad \forall \tilde{\omega}. \tag{9}$$

With the introduction of $\omega_\ell$, the variational formulation of the conservation law (2) becomes
$$-(\nabla \cdot (\ell\omega_\ell), \tilde{r})_\Omega + (S\sigma, \tilde{r})_\Omega = 0, \qquad \forall \tilde{r}. \tag{10}$$

Equations (7)-(10) now constitute the continuous, variational formulation in terms of $(\sigma, \omega_\ell, u, r_s)$. By inspecting the equations, we now determine the minimal continuity required of the four variables, which serves as a guideline for choosing appropriate finite element spaces.

First, we note that there are no differential operators acting on the rotation and displacement variables. These can therefore be represented by tuples of discontinuous, piecewise polynomials of order $k$, denoted by $\mathbb{P}_k$. Second, we observe that the Cauchy and couple stresses are subjected to the divergence operator, but not a full gradient. We will therefore only need to enforce continuity of their normal components across element boundaries. We achieve this by representing each row of the stress and couple stress tensors suitably chosen elements from either the Brezzi-Douglas-Marini ($\mathbb{BDM}_{k+1}$) or the Raviart-Thomas finite element space ($\mathbb{RT}_k$) of order $k$, cf. (Boffi, Brezzi and Fortin 2013) for more details on these spaces.

In the lowest order case, the finite element problem is as follows: find $(\sigma, \omega_\ell, u, r_s) \in \mathbb{BDM}_1^3 \times \mathbb{RT}_0^3 \times \mathbb{P}_0 \times \mathbb{P}_0$ such that
$$\begin{aligned}
(A\sigma, \tilde{\sigma})_\Omega + (u, \nabla \cdot \tilde{\sigma})_\Omega - (r_s, S\tilde{\sigma})_\Omega &= 0 \\
((2\mu)^{-1}\omega_\ell, \tilde{\omega})_\Omega + (r_s, \nabla \cdot (\ell\tilde{\omega}))_\Omega &= 0 \\
(-\nabla \cdot \sigma, \tilde{u})_\Omega &= (f, \tilde{u})_\Omega \\
-(\nabla \cdot (\ell\omega_\ell), \tilde{r})_\Omega + (S\sigma, \tilde{r})_\Omega &= 0
\end{aligned} \tag{11}$$
for all $(\tilde{\sigma}, \tilde{\omega}, \tilde{u}, \tilde{r}) \in \mathbb{BDM}_1^3 \times \mathbb{RT}_0^3 \times \mathbb{P}_0^3 \times \mathbb{P}_0^3$. The two-dimensional analogue of system (11) is posed on the space $\mathbb{BDM}_1^2 \times \mathbb{RT}_0 \times \mathbb{P}_0^2 \times \mathbb{P}_0$. These finite element spaces are well-defined for simplicial and Cartesian grids.

In the Cauchy limit, $\ell = 0$ eliminates the coupling terms between the second and last equations in (11), thereby forcing $\omega_\ell = 0$ strongly and $S\sigma = 0$ in a weak sense. The scheme then coincides with the mixed finite element method for elasticity with weakly imposed symmetry, proposed in (Arnold, Falk and Winther 2007). We note that the simpler space $\mathbb{RT}_0$ is not a suitable choice for the Cauchy stress since it would not be stable in the Cauchy limit. Stability in incompressible limit $\lambda \to \infty$, on the other hand, is captured as a consequence of the definition of the operator $A$ (Boffi, Brezzi and Fortin 2013).

The finite element discretization (11) was shown to be stable and convergent with optimal rates on shape-regular, simplicial grids in (Boon, Duran and Nordbotten 2024). We emphasize that, after negating the final two equations, the system is symmetric and therefore amenable to iterative solvers such as MinRes. However, it is a saddle-point system which is typically more challenging to solve compared to a positive definite system. Moreover, it is a relatively large system, containing 6 degrees of freedom per mesh element for the displacements and rotations (3 in 2D) and 12 degrees of freedom per face for the stresses (5 in 2D).

## 1.2 A two-point stress approximation finite volume method

The second discretization method we consider for the system (1)-(4) is the finite volume method proposed in (Nordbotten and Keilegavlen 2024). In this case, we rewrite the equations to a conservation form and subsequently approximate the differential operators using appropriate finite differences. Again, we pay particular attention to the Cauchy and incompressibility limits.

Similar to the mixed finite element method of Section 2.1, we construct a method that represents the primary variables such as displacements and rotations as cell-centred variables and defines the secondary variables, i.e. the Cauchy and couple stress, on the mesh faces. We emphasize that in such an approach for linearized elasticity, the first obstacle is typically the discretization of the symmetric gradient in (6). However, due to the reformulation in terms of the rotation and couple stress, we have conveniently avoided this difficulty.

Thus, the first obstacle we encounter concerns the term $\lambda(\nabla \cdot u)I$ in (2), which cannot directly be evaluated on mesh faces using the displacement degrees of freedom in the two adjacent cells. To handle this, we introduce two additional variables to the system, namely the *solid pressure* $p$ as a primary variable and the *displacement flux* $v$ as a secondary variable, given by

$$p := \lambda(\nabla \cdot u), \qquad v := u. \tag{12}$$

Additionally, we introduce the *rotation stress* $r$ as a primary variable, and the *total rotation* $\tau$ as a secondary variable:

$$r := 2\mu r_s, \qquad \tau := \frac{1}{2}\omega + S^*u. \tag{13}$$

We obtain the relation between the total stress and total rotation through the following calculation:

$$\nabla \cdot \tau = \frac{1}{2}\nabla \cdot \omega + \nabla \cdot S^*u = \frac{1}{2\mu}S(\sigma) + \nabla \cdot S^*u = S(\nabla u + S^*r_s) + \nabla \cdot S^*u = 2r_s = \mu^{-1}r \tag{14}$$

Here, we used the definition of $\tau$ in the first equality, followed by the conservation law (3), Hooke's law (2), and finally the identities $S(\nabla u) = -\nabla \cdot S^*u$ and $SS^*r_s = 2r_s$, which can be verified by direct computation. The same relation holds in 2D for scalar $r$ and vector-valued $\tau$, following a similar derivation.

We now describe the Cosserat system as three conservation equations and three constitutive laws in terms of the primary variables $(u, r, p)$ and secondary variables $(\sigma, \tau, v)$. First, the conservation equations are given by momentum balance (1), relation (14), and the definition of solid pressure $p$ in (12):

$$\nabla \cdot \begin{pmatrix} \sigma \\ \tau \\ v \end{pmatrix} - \begin{pmatrix} 0 & \square & \square \\ \square & \mu^{-1} & \square \\ \square & \square & \lambda^{-1} \end{pmatrix} \begin{pmatrix} u \\ r \\ p \end{pmatrix} = \begin{pmatrix} -f \\ 0 \\ 0 \end{pmatrix} \tag{15}$$

The discretization of (15) is now immediate by defining the primary variables on the mesh cells and the secondary variables on the mesh faces. In particular, we evaluate (15) by integrating over a control volume, i.e. mesh cell, $K_i$ and using the divergence theorem:

$$\sum_{k:\, \varsigma_k \subseteq \partial K_i} \Delta_{i,k} \begin{pmatrix} \sigma_k \\ \tau_k \\ v_k \end{pmatrix} - |K_i| \begin{pmatrix} 0 \\ \mu_i^{-1} r_i \\ \lambda_i^{-1} p_i \end{pmatrix} = |K_i| \begin{pmatrix} -f_i \\ 0 \\ 0 \end{pmatrix} \tag{16}$$

Here the incidence matrix $\Delta_{i,k} = \pm 1$ is positive if the unit normal $n_k$ associated with face $\varsigma_k$ is outward with respect to $K_i$, and negative otherwise. $|K_i|$ denotes the volume of cell $K_i$. Here, and in the following, the unknowns of any secondary variable describe the normal components on the faces, i.e. $\sigma_k = n_k \cdot \sigma$.

We continue with the constitutive laws that define the secondary in terms of the primary variables. They are Hooke's law (2), the definition of $\tau$ from (13) combined with (4), and the definition of $v$ from (12):

$$\begin{pmatrix} \sigma \\ \tau \\ \upsilon \end{pmatrix} = \begin{pmatrix} 2\mu\nabla & S^* & I \\ S^* & \ell^2\nabla & \square \\ I & \square & \square \end{pmatrix} \begin{pmatrix} u \\ r \\ p \end{pmatrix} \tag{17}$$

To evaluate (15) on the mesh faces, we introduce some geometric quantities. For a connected cell-face pair $(K_i, \varsigma_k)$, let $\delta_k^i \coloneqq (\Delta_{i,k} n_k) \cdot (x_k - x_i)$ denote the projected distance from the centre of cell $K_i$ to the centre of face $\varsigma_k$. Let its reciprocal be given by $\delta_k^{-i} \coloneqq \frac{1}{\delta_k^i}$. We then define the weighted averaging operators $\Xi$ and $\tilde{\Xi}$ as

$$\Xi_{k,i} \coloneqq \frac{\mu_i \delta_k^{-i}}{\sum_j \mu_j \delta_k^{-j}}, \qquad \tilde{\Xi}_{k,i} \coloneqq 1 - \Xi_{k,i}.$$

Here, the sum in the denominator is taken over all cells $K_j$ adjacent to face $\varsigma_k$. Letting these neighbouring cells be indexed $i$ and $j$ for notational brevity, we moreover define the face-orthogonal distance $\delta_k$ and the $\mu$-weighted distance $\delta_k^\mu$ as:

$$\delta_k \coloneqq \delta_k^i + \delta_k^j, \qquad \delta_k^\mu \coloneqq \frac{1}{2}(\mu_i \delta_k^{-i} + \mu_j \delta_k^{-j})^{-1}$$

Finally, we define the harmonic averages of the material parameters on the faces:

$$\bar{\mu}_k \coloneqq \delta_k \frac{\mu_i \delta_k^{-i} \mu_j \delta_k^{-j}}{\mu_i \delta_k^{-i} + \mu_j \delta_k^{-j}}, \qquad \overline{\ell^2}_k \coloneqq \delta_k \frac{\ell_i^2 \delta_k^{-i} \ell_j^2 \delta_k^{-j}}{\ell_i^2 \delta_k^{-i} + \ell_j^2 \delta_k^{-j}}$$

With these quantities defined, the discretization of (17) on a face $\varsigma_k$ becomes:

$$\begin{pmatrix} \sigma_k \\ \tau_k \\ \upsilon_k \end{pmatrix} = |\varsigma_k| \sum_{j:\, \varsigma_k \subseteq \partial K_j} \begin{pmatrix} -\delta_k^{-1} 2\bar{\mu}_k \Delta_{j,k} & -(S^* n_k)\tilde{\Xi}_{k,j} & n_k \tilde{\Xi}_{k,j} \\ -(S^* n_k)\Xi_{k,j} & -\delta_k^{-1} \overline{\ell^2}_k \Delta_{j,k} & \square \\ n_k \Xi_{k,j} & \square & -\delta_k^\mu \Delta_{j,k} \end{pmatrix} \begin{pmatrix} u_j \\ r_j \\ p_j \end{pmatrix} \tag{18}$$

For the derivation of (18) and details concerning boundary conditions, we refer to (Nordbotten and Keilegavlen 2024). Note that the terms containing $-\delta_k^{-1} \Delta_{j,k}$ represent a two-point finite difference that compares the values of the primary variables in adjacent cells and divides by the distance. The resulting stencil is minimal in the sense that it corresponds to a five-point stencil on three-dimensional simplicial grids, respectively a four-point stencil in 2D.

The two-point stress approximation finite volume method for (1)-(4) is now given by substituting (18) in (16). We emphasize that the incompressible limit $\lambda \to \infty$ and the Cauchy limit $\ell \to 0$ are naturally handled by this scheme as the relevant terms simply vanish in equations (16) and (18), respectively. Robustness of the method with respect to these limits is shown rigorously in (Nordbotten and Keilegavlen, Two-point stress approximation: A simple and robust finite volume method for linearized (poro-)mechanics and Stokes flow 2024).

Since the method only includes the primary variables, the resulting system consists of only 7 degrees of freedom per mesh cell (5 in 2D). While the method has been shown to be consistent and convergent under mild assumptions on the grid, it is applicable to any polytopal grid that satisfies $\delta_k > 0$.

### 1.3 A multi-point stress finite volume method
The multi-point stress approximation (MPSA) finite volume method was proposed about a decade ago (Nordbotten 2014), and has since been made available in multiple open-source simulation tools, such as the MRST (Lie 2019), PuMA (Ferguson, et al. 2021) and PorePy . The MPSA methods are known to be convergent on a broad class of grids, including simplexes as considered in the numerical experiments below (Nordbotten 2016). Of the different variations of the MPSA methods, we will here consider the so-called MPSA-W method, which is based on the techniques related to weak imposition of symmetry detailed in Section 2.1 (Keilegavlen and Nordbotten, Finite volume methods for elasticity with weak symmetry 2017).

The construction of the MPSA-W discretization stencil is somewhat technical, and we will not repeat it herein. However, we emphasize some of the qualitative properties of the discretization. In comparison to the MFEM and TPSA methods discussed in the previous sections, the MPSA-W method has a minimum number of degrees of freedom, consisting only of the vector of displacement (three degrees of freedom in 3D, and two in 2D) for each cell of the grid. However, this advantage is offset by a significantly larger stencil than the TPSA method, consisting of all neighbors of all corners of each cell. The exact size of the stencil depends on the number of simplexes meeting at each vertex, and is thus not unique. However, on logically Cartesian grids in 3D, this size of the stencil evaluates to containing 27 cells. Furthermore, the calculation of this discretization stencil requires the solution of a local problem for each corner of the grid, which itself will consist of a matrix problem on the order of 100 degrees of freedom, with multiple right-hand sides.

## 2. Numerical experiments

To compare the different numerical methods, we consider three numerical test cases. The governing equations for the first two cases are those of a linear elastic medium, where we consider robustness of the numerical methods under, respectively, the incompressible limit in a spatially homogeneous material, and the impact of material heterogeneities. In the third case, we consider the extended formulation of a Cosserat material, and consider a spatially varying Cosserat parameter $\ell$.

For the linearly elastic problems, we consider the mixed finite element formulation with three fields (as proposed in (Arnold, Falk and Winther 2007)) as well as TPSA and MPSA, while for the Cosserat problem, we consider the mixed finite element formulation with four fields (as proposed in (Boon, Duran og Nordbotten 2024)) as well as TPSA.

The computational meshes are for all examples simplicial tessellations of the unit square. However, to resolve the parameter heterogeneities, different meshes are considered in the three cases. This has the implication that the quality of the mesh deteriorates from the first to last example. Importantly, while essentially all triangles in example 3.1 are acute, this does not hold for example 3.3.

For each experiment, we impose a manufactured solution for the displacement $u = (u_x, u_y)$. The accuracy of the methods is measured in terms of the (cell center) displacement fields and the stress, with error norms computed as compared to the manufactured solution projected onto the grid. I.e. for $\pi_u$ being the $L^2$ projections onto the piecewise constants on cells and and $\pi_\sigma$ being the projection onto RT0, we measure the errors as:

$$e_u = \|u_h - \pi_u u\|_\Omega \qquad \text{and} \qquad e_\sigma = \|\sigma_h - \sigma\|_\omega$$

### 2.1 Linear elastic medium with homogeneous parameters

We first consider a purely elastic material with spatially homogenous parameters, where we set $\mu = 1$ and vary $\lambda \in \{10, 10^2, 10^4, 10^8\}$. This allows us to assess the robustness of the methods with respect to the physically important incompressible limit.

We report convergence of the numerical methods towards the expression

$$u = (u_x, u_y) = \left(\frac{\partial \psi}{\partial y}, -\frac{\partial \psi}{\partial x}\right), \quad \text{where} \quad \psi = \sin^2(2\pi x) \sin^2(2\pi y),$$

By construction the displacement $u$ is divergence free and zero on the boundaries of the domain.

The results for the MFEM, TPSA and MPSA methods are reported in Figure 1. As we can see, all three methods in general enjoy second-order convergence in displacement for this example, and first-order convergence in stress. We note that the second-order convergence in displacement is a kind of "super-convergence" phenomena, as in particular, only first-order convergence should be expected for TPSA (Nordbotten og Keilegavlen 2024). Both the MFEM and TPSA methods are furthermore fully robust in the incompressible limit, to the extent where the various curves are hardly distinguishable in the plot. On the other hand, the MPSA method is robust for nearly incompressible materials ($\lambda = 10^4$), but suffers slightly for the extreme case of $\lambda = 10^8$. A remedy for the MPSA method in the incompressible limit has been proposed previously (Keilegavlen and Nordbotten, Finite volume methods for elasticity with weak symmetry 2017), but is not considered herein.

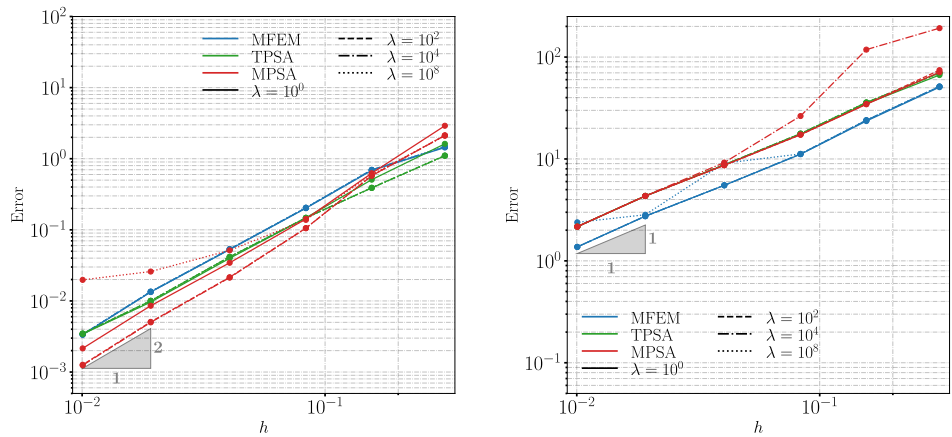

*Figure 1* Convergence history for the three methods considered Example 1. Left panel shows error data for $e_u$, right panel for $e_\sigma$. The colors differentiate the methods, and the line style differentiates the value of the Lamé parameter $\lambda$.

## 2.2 Linear elastic medium with parameter heterogeneities

To study the robustness of the methods in the presence of discontinuous material coefficients, we introduce a heterogeneity in the upper right corner of the simulation domain. We define the characteristic function defining the upper right corner of the unit square as

$$\chi_L = \begin{cases} 1 & \min(x,y) > 1/2 \\ 0 & \text{otherwise} \end{cases},$$

and set the Lame parameters as heterogeneous according to $\mu = (1 - \chi_L) + \chi_L \kappa$, and $\lambda = \mu$. The manufactured solution for the displacement is set to

$$u = \begin{pmatrix} \sin(2\pi x) y (1-y) \left(\frac{1}{2} - y\right) \\ \sin(2\pi y) x (1-x) \left(\frac{1}{2} - x\right) \end{pmatrix} / ((1 - \chi_L) + \chi_L \kappa)$$

The analytical displacement is constructed to have a sharp discontinuity in displacement gradient at the boundary of the material contrast, but still produce a continuous displacement and normal component of stress.

The mesh is constructed to conform to the material boundary.

The results for the MFEM, TPSA and MPSA methods are reported in Figure 2. As in Example 1, all three methods again preserve second-order convergence in displacement and first-order convergence in stress. This holds true independent of the presence of a material discontinuity of four orders of magnitude ($\kappa = 10^{\pm 4}$). We note that as the mesh is constructed to conform to the material boundary, there is a presence of obtuse angles in some of the triangles of a mesh. As expected from theory, this impacts the convergence rate of the displacement calculated by the TPSA method slightly (as can be seen in the left figure), but does not affect the two other methods.

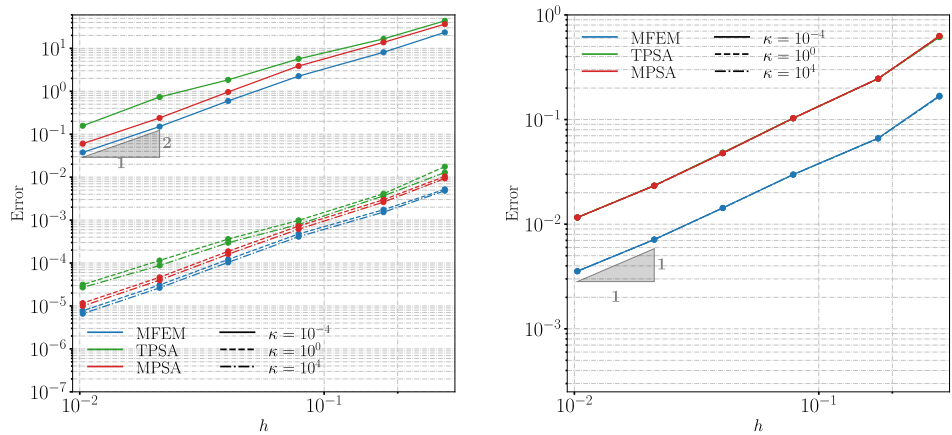

*Figure 2* Convergence history for the three methods considered Example 2. Left panel shows error data for $e_u$, right panel for $e_\sigma$. The colors differentiate the methods, and the line style differentiates the value of the material heterogeneity parameter κ. The data for the TPSA method in the right panel is hidden behind the MPSA data.

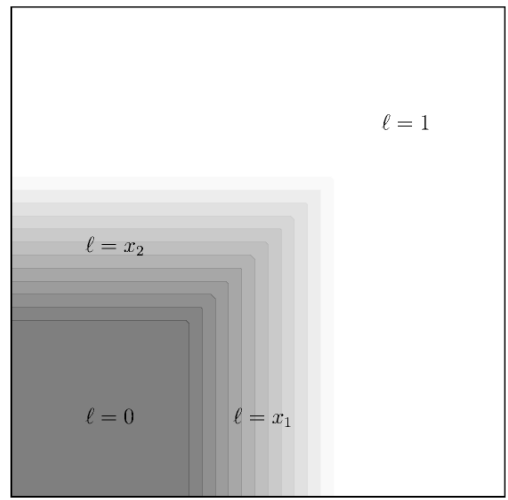

*Figure 3* Illustration of the spatially varying $\ell$ used in Example 3. Recall that $\ell = 0$ corresponds to a regular linearly elastic medium.

### 2.3 Composite Cosserat/elastic materials

Our final example considers a composite material, where parts of the material is modelled by linearized elasticity, while parts are modelled by the Cosserat description including the couple stress $\omega$. This is achieved by using a spatially varying length scale parameter $\ell$ that is non-zero in parts of the domain, according to

$$\ell = \min(1, \max(0, \max(3x - 1, 3y - 1)))$$

This is illustrated in Figure 3, and we emphasize that the mesh is adapted to resolve exactly the transition between the linear and constant regions of $\ell$. This leads to the presence of a higher number of obtuse angles in the resulting triangulation.

For this example, the Lamé parameters are both set to unity. The manufactured solution for displacement is chosen as
$$u = (u_x, u_y) = (\sin(2\pi x) y (1 - y), \sin(2\pi y) x (1 - x))$$

To ensure a non-trivial couple stress in the Cosserat region, we define the rotation field as:
$$r = xy(1 - x)(1 - y)$$

The results for the MFEM and TPSA methods are reported in Figure 4 (there is to date not published a generalization of MPSA to Cosserat materials). As in Examples 1 and 2, the MFEM method preserves the second-order convergence of displacement and first-order convergence of stress. The TPSA method remains convergent for this problem, although the convergence rates are reduced. While it cannot be deduced directly from Figure 4, we have verified independently that this reduction in convergence rates is due to the lack of an acute triangulation, not the material parameters.

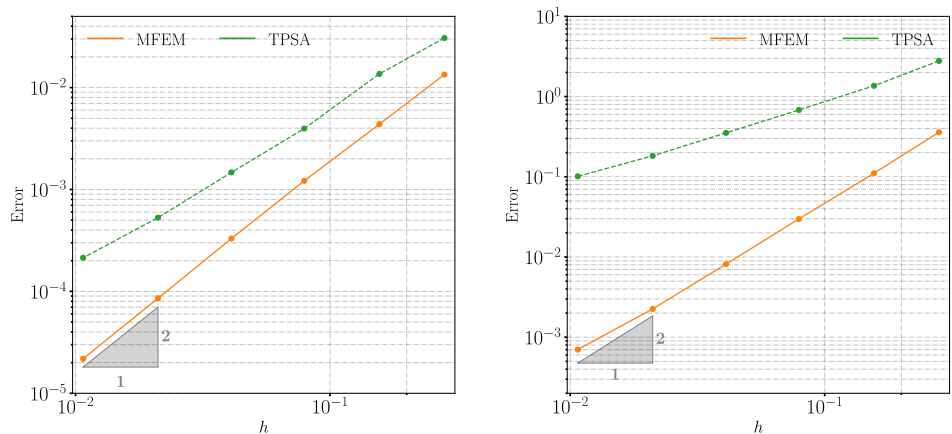

*Figure 4* Convergence history for the two methods considered Example 3. Left panel shows error data for $e_u$, right panel for $e_\sigma$.

### 2.4 Grid size and degrees of freedom

The comparisons in Sections 3.1-3.3 are based on computed L2 error relative to grid size. However, the grid size is only one of several aspects that may impact the computational cost, which ultimately will also depend on implementation and hardware. Equally relevant may be the size of the solution vector, determined by the number of degrees of freedom (DOF).

Recall that the MPSA method has two DOFs per cell (displacement), while TPSA has four DOFs per cell (displacement, rotation and solid pressure). The number of DOFs for MFEM depends on the problem type. For the three-field formulation approximating the equations of elasticity, the MFEM has on average nine DOFs per cell (displacement, rotation, and four DOFs of stress per face). For the four-field formulation approximating the equations of Cosserat materials, MFEM has on average 10.5 DOFs per cell (an additional one DOFs of couple stress per face).

Two other aspects impacting the performance of the discretization are the assembly time of the linear system, matrix fill-in and the condition number of the resulting system. We have not studied these issues systematically herein, but make the following comments.

In terms of matrix assembly, both the MFEM and TPSA methods are relatively cheap, as they only involve direct evaluation of integrals and simple algebraic expressions. In contrast the MPSA method requires a somewhat expensive assembly of the linear system, essentially involving a local static condensation of stress degrees of freedom on faces.

The linear systems for MPSA, TPSA and MFEM have very different structure and characteristics, which must be treated by different (iterative) solver strategies. As such, it is not at present possible to give a simple assessment of which systems are most efficient in applications.

3.  Conclusions

We have compared the approximation properties three recently developed methods for approximating the solution to linearly elastic problems. All three methods share the attractive characteristic that they explicitly balance momentum, which may be attractive for problems where surface stress (traction) is important, such as for fractures, faults and in implementing general boundary conditions.

The MFEM method is in a certain sense the most complex method, requiring a high number of degrees of freedom and resulting in a saddle-point type system. As presented in the cited references, the method is also limited to rather simple grid types (simplexes and nice quadrilaterals). On the other hand, within its limitations, it is robust for all parameter regimes and heterogeneities considered herein.

The TPSA method is conceptually the simplest method, using only two cells to obtain an explicit approximation to the stress across any face of the grid. This leads to a method that is simple and fast to implement, and applicable to very general grid types. It is also robust for all parameter regimes and heterogeneities considered herein. However, as with all two-point methods, the convergence will depend on the shape of the cells forming the grid.

Finally, the MPSA can be considered as a compromise between the TPSA and MFEM methods. It is related to TPSA in the sense of having explicit and local approximations to the stress across any face of the grid (albeit involving more than two neighboring cells). The MPSA method is also related to a certain MFEM formulation through various quadrature rules (Ambartsumyan, et al. 2020). This compromise enhances the robustness of MPSA over TPSA with respect to unfavorable cell shapes, but comes at a penalty in terms of the complexity of assembling the method and the matrix fill-in of the resulting linear system.

**Acknowledgements**

OD was supported by the European Research Council (ERC), under the European Union's Horizon 2020 research and innovation program (grant agreement No 101002507). JMN was supported in part through Norwegian Research Center (NORCE) and the project Expansion of Resources for CO2 Storage on the Horda Platform (ExpReCCS) grant nr. 336294 financed through the Research Council of Norway, Equinor, Norske Shell and Wintershall DEA Norge.